\documentclass[11pt]{article}
\usepackage{eurosym}
\usepackage{amsthm,amssymb,amsmath}
\usepackage{epsfig}
\usepackage{color}
\usepackage{pgf,tikz}

\newtheorem{theorem}{Theorem}
\newtheorem{corollary}[theorem]{Corollary}

\newtheorem{obs}[theorem]{Observation}
\newtheorem{proposition}[theorem]{Proposition}

\theoremstyle{definition}

\theoremstyle{remark}

\DeclareMathOperator{\opt}{\emph{Opt}}

\begin{document}

\title{Maximal double Roman domination in graphs}
\date{}
\author{H. Abdollahzadeh Ahangar$^{(1)}$, M. Chellali$^{(2)}$, \and S.M.
Sheikholeslami$^{(3)}$ and J.C. Valenzuela-Tripodoro$^{(4)}$\vspace{2mm} \\
$^{(1)}${\small Department of Mathematics}\\
{\small \ Babol Noshirvani University of Technology}\\
{\small \ Shariati Ave., Babol, I.R. Iran, Post Code:47148-71167.}\\
{\small \ ha.ahangar@nit.ac.ir\vspace{2mm}}\\
$^{(2)}${\small LAMDA-RO Laboratory, Department of Mathematics}\\
{\small \ University of Blida}\\
{\small \ B.P. 270, Blida, Algeria}\\
{\small \ m\_chellali@yahoo.com\vspace{2mm}}\\
$^{(3)}${\small Department of Mathematics }\\
{\small \ Azarbaijan Shahid Madani University }\\
{\small \ Tabriz, I.R. Iran}\\
{\small \ s.m.sheikholeslami@azaruniv.edu\vspace{3mm}}\\
$^{(4)}${\small Department of Mathematics, University of C\'{a}diz, Spain.}
\\
\ {\small jcarlos.valenzuela@uca.es\vspace{5mm}}}
\maketitle

\begin{abstract}
A maximal double Roman dominating function (MDRDF) on a graph $G=(V,E)$ is a
function $f:V(G)\rightarrow \{0,1,2,3\}$ such that \textrm{(i) }every vertex
$v$ with $f(v)=0$ is adjacent to least two vertices { assigned $2$
or to at least one vertex assigned $3,$} \textrm{(ii) }every vertex $v$ with
$f(v)=1$ is adjacent to at least one { vertex assigned $2$ or $3$}
and \textrm{(iii) }the set $\{w\in V|~f(w)=0\}$ is not a dominating set of $G
$. The weight of a MDRDF is the sum of its function values over all
vertices, and the maximal double Roman domination number $\gamma _{dR}^{m}(G)
$ is the minimum weight of an MDRDF on $G$. {In this paper, we initiate the
study of maximal double Roman domination. We first show that the problem of
determining }$\gamma _{dR}^{m}(G)$ {is NP-complete for bipartite, chordal
and planar graphs. But it is solvable in linear time for bounded
clique-width graphs including trees, cographs and distance-hereditary
graphs. Moreover, we establish various relationships relating }$\gamma
_{dR}^{m}(G)$ to some domination parameters. {For the class of trees, we
show that for every tree }$T$ {of order }$n\geq 4,$ $\gamma _{dR}^{m}(T)\leq
\frac{5}{4}n$ {and we characterize all trees attaining the bound. Finally,
the exact values of }$\gamma _{dR}^{m}(G) $ {are given for paths and cycles.
}\newline

\end{abstract}

\section{Introduction}

In this paper, $G$ is a simple graph with vertex set $V=V(G)$ and edge set $%
E=E(G)$. The order $|V|$ of $G$ is denoted by $n$. For every vertex $v\in V$%
, the \emph{open neighborhood} $N(v)$ is the set $\{u\in V(G):uv\in E(G)\}$
and the \emph{closed neighborhood} of $v$ is the set $N[v]=N(v)\cup \{v\}$.
The \emph{degree} of a vertex $v\in V$ is $\deg (v)=|N(v)|$. The \emph{%
minimum} and \emph{maximum degree} of a graph $G$ are denoted by $\delta
=\delta (G)$ and $\Delta =\Delta (G)$, respectively. A \textit{leaf} of $G$
is a vertex of degree one, while a \textit{support vertex }of $G$ is a vertex adjacent
to a leaf.

A set $S\subseteq V$ in a graph $G$ is called a \textit{dominating set} if
every vertex of $G$ is either in $S$ or adjacent to a vertex of $S.$ The
\textit{domination number} $\gamma (G)$ equals the minimum cardinality of a
dominating set in $G$. A dominating set $D$ is said to be a \textit{maximal
dominating set} (MDS) if $V-D$ is not a dominating set of $G$. The \emph{%
maximal domination number }$\gamma _{m}(G)$ is the minimum cardinality of an
MDS of $G$. {Maximal domination was first defined by }Kulli and Janakiram
\cite{kj} in 1997.

For a graph $G$ and a positive integer $k$, let $f:V(G)\rightarrow
\{0,1,2,\dots ,k\}$ be a function, and let $(V_{0},V_{1},V_{2},\ldots
,V_{k}) $ be the ordered partition of $V=V(G)$ induced by $f$, where $%
V_{i}=\{v\in V:f(v)=i\}$ for $i\in \{0,1,\ldots ,k\}$. {Since} $f $ {is
determined by these sets, we will write }$f=(V_{0},V_{1},...,V_{k})$ (or $%
f=(V_{0}^{f},V_{1}^{f},...,V_{k}^{f})$ to refer to $f$). {Moreover, the
weight of }$f$ {is given by }$f(V(G))=\sum_{u\in V(G)}f(u).$

A function $f=(V_{0},V_{1},V_{2})$ is a \textit{Roman dominating function}
(RDF) on $G$ if every vertex $u\in V_{0}$ {has at least one neighbor in }$%
V_{2}.$ The \textit{Roman domination number} $\gamma _{R}(G)$ is the minimum
weight of an RDF on $G$. Roman domination was introduced {in 2004} by
Cockayne et al. \cite{CDHH} {and since then }more than 200 papers have been
published on this topic, where several new variations were introduced. {In
particular, those which interest us in this paper: }maximal Roman domination
\cite{asstv} and double Roman domination \cite{bhh}. For more details on
Roman domination and its variants we refer the reader to the recent three papers \cite{amc1, amc2, acsv}, two
book chapters \cite{Ch1, Ch2} and surveys papers \cite{Ch3}.

{An RDF function }$f=(V_{0},V_{1},V_{2})$ {is a }\textit{maximal Roman
dominating function} (MRDF) {}on $G$ {if }$V_{0}$ is not a dominating set of
$G$. The \emph{maximal Roman domination number} {(\emph{maximal RD-number},
for short) }$\gamma _{mR}(G)$ of $G$ equals the minimum weight of an MRDF of{%
\ }$G$. {Maximal Roman dominating functions were introduced in \cite{asstv}
motivated by maximal dominating sets}{ , for more see \cite{amds}}.

As defined in \cite{bhh}, a function $f=(V_{0},V_{1},V_{2},V_{3})$ is a
\textit{double Roman dominating function} (DRDF) on a graph $G$ if\textrm{\ }%
the following two conditions hold: \textrm{(i) }{every vertex in }$V_{0}$ {{%
 must have a neighbor in }$V_{3}$ or at least two neighbors in $%
V_{2}$; \textrm{(ii) }{every vertex in }$V_{1}$ {must have a neighbor }in $%
V_{2}\cup V_{3}.${\ }The \textit{double Roman domination number }{(\emph{%
DRD-number}, for short) }$\gamma _{dR}(G)$ equals the minimum weight of a
DRDF on $G$. }

Let us recall that the main motivation for introducing double Roman
domination by Beeler et al. \cite{bhh} was to strengthen the defense of the
Roman Empire where three legions can be deployed at a given location, and
thus offering a high level of defense ensuring that any attack can be
defended by at least two legions. {What we propose in this paper is a
stronger version of double Roman domination in the sense that we adopt this
strategy of defense but we also make sure that there exists a location that
keeps at least one legion {\ that must not defend any defenceless neighbor. }%
}

\begin{figure}[tbp]
\center
\begin{tikzpicture}
\node [draw, shape=circle,fill=black,scale=0.5] (b1) at  (2,5) {};
\node [draw, shape=circle,fill=black,scale=0.5] (b2) at  (3,5) {};
\node [draw, shape=circle,fill=black,scale=0.5] (b3) at  (4,5) {};
\node [draw, shape=circle,fill=black,scale=0.5] (b4) at  (5,5) {};
\node [draw, shape=circle,fill=black,scale=0.5] (b5) at  (7,5) {};
\node [draw, shape=circle,fill=black,scale=0.5] (b6) at  (8,5) {};
\node [draw, shape=circle,fill=black,scale=0.5] (b7) at  (9,5) {};
\node [draw, shape=circle,fill=black,scale=0.5] (b8) at  (10,6) {};
\node [draw, shape=circle,fill=black,scale=0.5] (b9) at  (10,8) {};
\node [draw, shape=circle,fill=black,scale=0.5] (a2) at  (8,7) {};
\node [draw, shape=circle,fill=black,scale=0.5] (a1) at  (3.5,7) {};
\node [draw, shape=circle,fill=black,scale=0.5] (c1) at  (3.5,3) {};
\node [draw, shape=circle,fill=black,scale=0.5] (c2) at  (2.5,1) {};
\node [draw, shape=circle,fill=black,scale=0.5] (c3) at  (4.5,1) {};

\draw(a1)--(a2);
\draw(a1)--(b1);\draw(a1)--(b2);\draw(a1)--(b3);\draw(a1)--(b4);
\draw(a2)--(b5);\draw(a2)--(b6);\draw(a2)--(b7);
\draw(a2)--(b8)--(b9)--(a2);
\draw(c1)--(b1);\draw(c1)--(b2);\draw(c1)--(b3);\draw(c1)--(b4);
\draw(c1)--(c2)--(c3)--(c1);
\node at (8,7.35)  {$3$} {};
\node at (3.5,7.35)  {$2$} {};
\node at (10.25,6)  {$0$} {};
\node at (10.25,8)  {$0$} {};
\node at (1.8,4.7)  {$0$} {};
\node at (2.8,4.7)  {$0$} {};
\node at (4.2,4.7)  {$0$} {};
\node at (5.2,4.7)  {$0$} {};
\node at (3,3)  {$2$} {};
\node at (2.25,1)  {$1$} {};
\node at (4.75,1)  {$1$} {};
\node at (7,4.7)  {$0$} {};
\node at (8,4.7)  {$0$} {};
\node at (9,4.7)  {$0$} {};

\end{tikzpicture}
\caption{Sink $(0)$; Reserve station$(1)$; Supply substation $(2)$ and
Supply station $(3)$.}
\label{electnetw}
\end{figure}
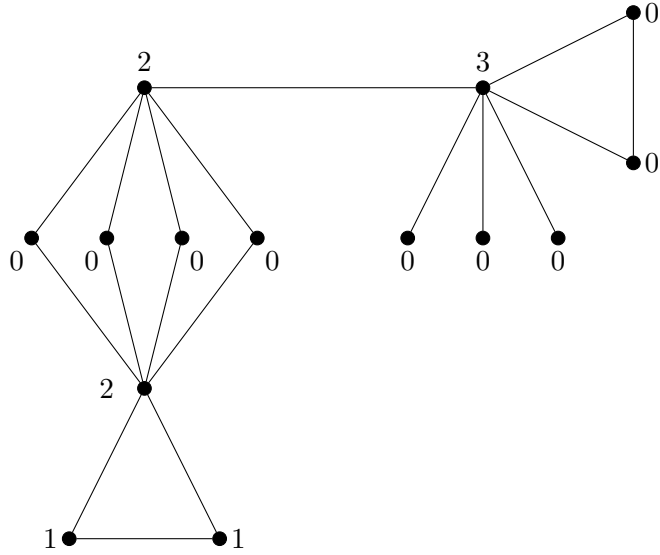

From a practical point of view, we may also apply this graph domination
strategy to the optimum design of connected networks. Some electrical
networks consist in four types of different interconnected elements, namely
sinks, reserve stations, supply substations and supply stations. Sinks are
elements that need to be connected either to a powerful supply station or to
two, less powerful, supply substations. Reserve stations must be connected
to a supply element and there must be at least one reserve station that is
not connected (i.e. not supplying energy) to any sink because it must
function as electricity storage (see Figure~\ref{electnetw} for an example).

{ Formally, we say that a double Roman dominating function}  $%
s=(V_{0},V_{1},V_{2},V_{3})$ on $G$ is a \emph{maximal double Roman
dominating function } (MDRDF) if $V_{0}$ is not a dominating set of $G$. The
\emph{maximal double Roman domination number} (\emph{maximal DRD-number},
for short) $\gamma _{dR}^{m}(G)$ equals the minimum weight of an MDRDF of $G$%
. A $\gamma _{dR}^{m}(G)$-\emph{function} is an MDRDF {of} $G$ with weight $%
\gamma _{dR}^{m}(G)$. Note that $\gamma _{R}(G)$-\emph{functions, }$\gamma
_{mR}(G)$-\emph{functions }and\emph{\ }$\gamma _{dR}(G)$-\emph{functions}
are similarly defined.

{In this paper, we initiate the study of maximal double Roman domination. A
DRDF }$f=(V_{0},V_{1},V_{2},V_{3})$ on $G$ is a {\emph{maximal double Roman
dominating function }(MDRDF) if }$V_{0}$ is not a dominating set of $G$. The
\emph{maximal double Roman domination number} {(\emph{maximal DRD-number},
for short) }$\gamma _{dR}^{m}(G)$ equals the minimum weight of an MDRDF of $G
$. A $\gamma _{dR}^{m}(G)$-\emph{function} is an MDRDF {of} $G$ with weight $%
\gamma _{dR}^{m}(G)$. {Note that }$\gamma _{R}(G)$-\emph{functions, }$\gamma
_{mR}(G)$-\emph{functions }and\emph{\ }$\gamma _{dR}(G)$-\emph{functions }{%
are similarly defined. }

{In this paper, we first show that the problem of determining }$\gamma
_{dR}^{m}(G)$ {is NP-complete for bipartite, chordal and planar graphs. Then
we show that this problem is solvable in linear time for bounded
clique-width graphs which include trees, cographs and distance-hereditary
graphs. Moreover, we establish various relationships relating the maximal
DRD-number of a graph to the DRD-number, maximal RD-number and the maximal
domination number. We also show that for trees of order }$n\geq 4,$ $\gamma
_{dR}^{m}(T)\leq \frac{5}{4}n$ {and a characterization of all trees
attaining the bound is provided. Finally, the exact values of the maximal
DRD-number for paths and cycles are established. }

{We close this section by the following result that gives two }properties of
{MDRDFs of a connected graph. { Recall that a vertex $x$ 
is said to be a private neighbor of a vertex $y\in D$ 
with respect to the set $D$  if $y$ is the only neighbor of $x$ in $D,$ that is $N(x)\cap
D=\left\{ y\right\} .$}

\begin{proposition}
{ Let $G$ be a connected graph.} Then:

\begin{itemize}
\item[(i)] For any $\gamma _{dR}^{m}$-function $(V_{0},V_{1},V_{2},V_{3})$
with $V_{3}\neq \emptyset ,$ every vertex of $V_{3}$ has a private
neighborhood in $V_{0}$ with respect to $V_{2}\cup V_{3}$.

\item[(ii)] { For any }$\gamma _{dR}^{m}$-function $%
(V_{0},V_{1},V_{2},V_{3})$ { with }$V_{0}\neq \emptyset ,$ $V_{0}$
does not dominate all $V_{1}\cup V_{2}$.
\end{itemize}
\end{proposition}

\textbf{Proof.} \textrm{(i)-} Let $f=(V_{0},V_{1},V_{2},V_{3})$ be a $\gamma
_{dR}^{m}$-function of $G$ and let $v\in V_{3}.$ If $v$ has no private
neighbor in $V_{0}$ {with respect to} $V_{2}\cup V_{3},$ then reassigning $v$
the value $2$ provides an MDRDF of $G$ {of weight }$\gamma _{dR}^{m}(G)-1$,
a contradiction. {Hence we get }\textrm{(i).}

\textrm{(ii)-} Let $(V_{0},V_{1},V_{2},V_{3})$ be a $\gamma _{dR}^{m}$%
-function of $G$ { with $V_{0}\neq \emptyset .$ 
Clearly, $V_{2}\cup V_{3}\neq \emptyset $ (since $V_{0}\neq
\emptyset $). If $V_{0}$ dominates all $V_{1}\cup V_{2},$ then
since $f$ is an MDRDF of $G,$ some vertex $x$ 
of $V_{3}$ must have no neighbor in $V_{0},$ 
contradicting item (i). } $\Box $


\section{Complexity results}

Our aim in this section is to study the complexity of the following decision
problem, to which we shall refer as MAXIMAL DOUBLE ROM-DOM:

\bigskip

\textbf{MAXIMAL\ DOUBLE ROM-DOM}

\textbf{Instance}: Graph $G=(V,E)$, positive integer $k$ ($k\leq2\left\vert
V(G)\right\vert $).

\textbf{Question}: Does $G$ have a maximal double Roman function of weight
at most $k$?

\bigskip

{First, to }show that this decision problem for maximal double Roman
domination is NP-complete, we use a polynomial time reduction from the
double Roman domination problem shown to be NP-complete for bipartite and
chordal graphs in \cite{acs1}, and for planar graphs in \cite{pj3}. \newline
\medskip

\textbf{DOUBLE\ ROMAN DOMINATION PROBLEM}

\textbf{INSTANCE:} A graph $G$ and a positive integer $k$ ($k\leq2\left\vert
V(G)\right\vert $).

\textbf{QUESTION:} Is $\gamma_{dR}^m(G)\leq k$?\newline

\begin{figure}[tbp]
\center
\begin{tikzpicture}
\node [draw, shape=circle,fill=black,scale=0.5] (e1) at  (5,-.7) {};
\node [draw, shape=circle,fill=black,scale=0.5] (a1) at  (5,0.5) {};
\node [draw, shape=circle,fill=black,scale=0.5] (a2) at  (5,1.5) {};
\node [draw, shape=circle,fill=black,scale=0.5] (a3) at  (5,2.2) {};
\node [draw, shape=circle,fill=black,scale=0.5] (a4) at  (4.5,2.2) {};
\node [draw, shape=circle,fill=black,scale=0.5] (a5) at  (5.5,2.2) {};
\draw(e1)--(a3);
\draw(a5)--(a2)--(a4);
\node at (4.6,.3)  {$a_1$} {};
\node at (4.6,1.3)  {$y_1$} {};
\node [draw, shape=circle,fill=black,scale=0.5] (e4) at  (6,-1.3) {};
\node [draw, shape=circle,fill=black,scale=0.5] (b1) at  (8,0) {};
\node [draw, shape=circle,fill=black,scale=0.5] (b2) at  (8.5,0.5) {};
\node [draw, shape=circle,fill=black,scale=0.5] (b3) at  (9,1.2) {};
\node [draw, shape=circle,fill=black,scale=0.5] (b8) at  (8.,1.2) {};
\node [draw, shape=circle,fill=black,scale=0.5] (b9) at  (8.5,1.2) {};
\draw(e4)--(b1)--(b2)--(b3)--(b2)--(b8)--(b2)--(b9);
\node at (8.1,-.3)  {$a_2$} {};
\node at (8.6,0.2)  {$y_2$} {};
\node [draw, shape=circle,fill=black,scale=0.5] (e3) at  (4,-1.3) {};
\node [draw, shape=circle,fill=black,scale=0.5] (c1) at  (2,0) {};
\node [draw, shape=circle,fill=black,scale=0.5] (c5) at  (1.5,.5) {};
\node [draw, shape=circle,fill=black,scale=0.5] (c6) at  (1.5,1.2) {};
\node [draw, shape=circle,fill=black,scale=0.5] (c7) at  (1.,1.2) {};
\node [draw, shape=circle,fill=black,scale=0.5] (c10) at  (2,1.2) {};
\draw(e3)--(c1)--(c5)--(c6)--(c5)--(c7)--(c5)--(c10);
\node at (1.8,-.2)  {$a_n$} {};
\node at (1.3,.3)  {$y_n$} {};
\node [draw, shape=circle,fill=black,scale=0.5] (e2) at  (5,-3.3) {};
\node [draw, shape=circle,fill=black,scale=0.5] (d1) at  (5,-4.6) {};
\node [draw, shape=circle,fill=black,scale=0.5] (d2) at  (5,-5.6) {};
\node [draw, shape=circle,fill=black,scale=0.5] (d3) at  (5,-6.3) {};
\node [draw, shape=circle,fill=black,scale=0.5] (d4) at  (4.5,-6.3) {};
\node [draw, shape=circle,fill=black,scale=0.5] (d5) at  (5.5,-6.3) {};
\node at (5.4,-4.5)  {$a_i$} {};
\node at (5.4,-5.5)  {$y_i$} {};
\draw(e2)--(d1)--(d2)--(d3)--(d2)--(d4)--(d2)--(d5);
\node [draw, shape=circle,fill=black,scale=0.3] (f1) at  (4.5,-2.3) {};
\node [draw, shape=circle,fill=black,scale=0.3] (f1) at  (4.3,-2.) {};
\node [draw, shape=circle,fill=black,scale=0.3] (f1) at  (4.7,-2.6) {};
\node [draw, shape=circle,fill=black,scale=0.3] (f1) at  (5.5,-2.3) {};
\node [draw, shape=circle,fill=black,scale=0.3] (f1) at  (5.7,-2.) {};
\node [draw, shape=circle,fill=black,scale=0.3] (f1) at  (5.3,-2.6) {};
\draw (5,-2) circle (1.8cm);
\node at (5, -1)  {$v_1$} {};
\node at (6.2, -1.6)  {$v_2$} {};
\node at (3.8, -1.6)  {$v_n$} {};
\node at (4.6, -3.3)  {$v_i$} {};
\end{tikzpicture}
\caption{NP-Completeness for bipartite, chordal and planar graphs.}
\label{Fig.1}
\end{figure}
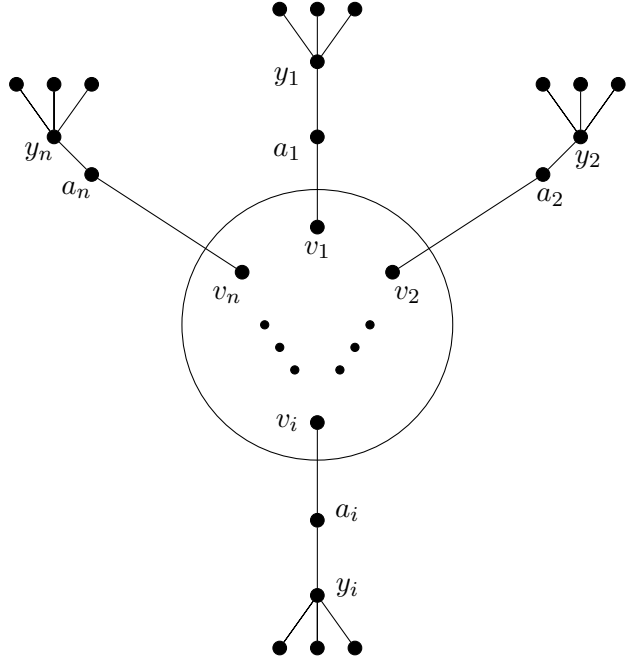

\begin{theorem}
Problem MAXIMAL\ DOUBLE ROM-DOM is NP-Complete for bipartite, chordal and
planar graphs\textbf{.}
\end{theorem}

\textbf{Proof. }Clearly, MAXIMAL\ DOUBLE ROM-DOM is a member of $\mathcal{NP}
$ since for a given function $f=(V_{0},V_{1},V_{2},V_{3})$ of a graph $G$ we
can check in polynomial time that $f$ is a DRDF for $G$ with weight at most $%
k$ and that $V_{0}$ does not dominate $G$.

Given a positive integer $k$ and a graph $G$ of order $n$, we construct a
graph $G^{\ast }$ by adding for each vertex $x_{i}$ a star $K_{1,4},$ with
center vertex $y_{i}$ and leaves $a_{i},b_{i},c_{i},d_{i},$ attached by $%
a_{i}x_{i}$ at $x_{i}.$ It is worth mentioning that $\left\vert V(G^{\ast
})\right\vert =6\left\vert V(G)\right\vert $ and $\left\vert E(G^{\ast
})\right\vert =\left\vert E(G)\right\vert +5\left\vert V(G)\right\vert ,$
and so $G^{\ast }$ can be constructed from $G$ in polynomial time. Morever,
it is clear that if $G$ is a bipartite, chordal or planar graph, then $%
G^{\ast }$ is also bipartite, chordal or planar, respectively.

Next, we shall show that $G$ has a DRDF $f$ with $\omega (f)\leq k$ if and
only if $G^{\ast }$ has an MDRDF $g$ with $\omega (g)\leq k+3n+1$. Suppose
that $f$ is a DRDF of $G$ with $\omega (f)\leq k.$ Define the function $%
g=(V_{0},V_{1},V_{2},V_{3})$ on $V(G^{\ast })$ by $g(x)=f(x)$ for all $x\in
V(G),$ and for every $i\in \{1,...,n\}$ (except for vertex $b_{1}$), let $%
g(b_{1})=1,$ $g(y_{i})=3,$ $g(a_{i})=g(b_{i})=g(c_{i})=g(d_{i})=0.$ Clearly,
$\omega (g)=\omega (f)+3n+1$ and thus $\omega (g)\leq k+3n+1.$ Moreover,
since vertex $b_{1}$ is adjacent to no { vertex assigned $0$}
under $g,$ we deduce that $g$ is an MDRDF of $G^{\ast }.$

Conversely, let $g$ be an MDRDF of $G^{\ast }$ with $\omega (g)\leq k+3n+1$.
We first note that for every $i,$ we have $%
g(y_{i})+g(b_{i})+g(c_{i})+g(d_{i})\geq 3$, and clearly the inequality is
strict if one of $b_{i},c_{i},d_{i}$ is assigned a positive value. Further
we may assume, without loss of generality, that $g(y_{i})=3$ and at least
one of $b_{i},c_{i},d_{i},$ say $b_{i},$ { is assigned $0$} under
$g$ for every $i.$ Suppose that $g(a_{j})\neq 0$ for some $j.$ If $%
g(a_{j})=1,$ then reassigning $a_{j}$ and $b_{j}$ the values $0$ and $1,$
respectively, provides another MDRDF of $G^{\ast }$ with the same weight $%
\omega (g).$ If $g(a_{j})=2,$ then the neighbor of $a_{j}$ belonging to $%
V(G),$ say $x_{j},$ is adjacent to at least one vertex in $V(G)$ with weight
$2$ or $3.$ Then reassigning $a_{j},b_{j},x_{j}$ the values $0,1$ and $1,$
respectively, provides another MDRDF of $G^{\ast }$ with the same weight $%
\omega (g).$ If $g(a_{j})=3,$ then it suffices to reassign $a_{j},b_{j},x_{j}
$ the values $0,1 $ and $2$ to obtain an MDRDF of $G^{\ast }$ with weight $%
\omega (g).$ Therefore, we can assume that $g(a_{i})=0$ for all $i.$
Clearly, the function $g$ restricted to $V(G)$ defines a DRDF $f$ on $G.$
Now since all $a_{i}$'s dominate $V(G)$ and $g$ is an MDRDF of $G^{\ast },$
we deduce that a leaf neighbor of some $y_{i},$ say $c_{i},$ is assigned a
positive value under $g$. It follows that $f$ is a DRDF $f$ on $G$ with $%
\omega (f)\leq \omega (g)-3n-g(c_{i})\leq k.$ $\ \Box $

\bigskip

{Despite the result that we have just proven, it is possible to show that
MAXIMAL\ DOUBLE ROM-DOM is solvable in linear time for a wide class of
graphs. Specifically, in what follows we prove that this problem is solvable
in linear time for all graphs with bounded clique-width (including trees,
cographs and distance-hereditary graphs for which the clique-width is
bounded by $3$). }

We need to make use of several concepts related to describing finite graphs
as logical structures. A $k$\textit{-presentation} or a $k$\textit{%
-expression} of a graph $G$ with vertex labels $\{1,2,\ldots ,k\}$ is a
finite set of some of the following operations

\begin{tabular}{rl}
$\bullet i(x): $ & To create a new vertex, $x,$ with an $i\in \{1,\ldots
,k\} $ assigned as a label. \\
$G_1\oplus G_2: $ & To create a new graph as the disjoint union of $G_1$ and
$G_2.$ \\
$\eta_{ij}(G): $ & To create all edges in $G$ that join $i$-vertices with $j$%
-vertices. \\
$\rho_{i\rightarrow j}(G):$ & To change the label of all $i$-vertices into
label $j.$%
\end{tabular}

The clique-width of a graph is defined as the minimum $k$ wich is needed to
define the graph $G$ by means of a $k$-expression. As an example, we can
make the bipartite graph $K_{2,3}$ from the following $2$-expression,
\begin{equation*}
\eta_{12} \left( \left( \left( \left( \bullet 1(a_1) \oplus \bullet 1(a_2)
\right) \oplus \bullet 2(b_1) \right) \oplus \bullet 2(b_2) \right) \oplus
\bullet 2(b_3) \right)
\end{equation*}
where $\{a_1,a_2\} \cup \{b_1,b_2,b_3\}$ are the partite set of vertices.

Moreover, MSOL($\tau_1$) represents the monadic second order logic with
quantification over subsets of elements of the logic structure $G(\tau_1)$.
Namely, $G(\tau_1)$ is the logic structure $<V(G),R>$ where $R$ is a binary
relation such that $R(x,y)$ is satisfied if and only if $x,y$ are adjacent
elements of $V(G).$

The class of $LinEMSOL(\tau)$ optimization problems includes those that can
be described as follows (see \cite{LKLP} for more details about the
restriction of the definition given in \cite{CMR} to finite graphs)
\begin{equation*}
\opt \; \left\{ \sum_{1\le i \le l} a_i |X_i| \;:\;
<G(\tau_1),X_1,\ldots,X_l > \; \vDash \theta(X_1,\ldots,X_l) \right\}
\end{equation*}
where $\theta$ is an MSOL($\tau_1$) formula that contains free set-variables
$X_1,\ldots,X_l, $ integers $a_i$ and \emph{Opt} is either $\min$ or $\max.$
See \cite{CMR, LKLP} for more details and formal definitions.

To our purpose, we make extensive use of a result regarding LinEMSOL
optimization problems by Courcelle et al. \cite{CMR}.

\begin{theorem}[Courcelle et al. \protect\cite{CMR}]
\label{courcelle}Let $k\in \mathbb{N}$ and let $\mathcal{C}$ be a class of
graphs of clique-width at most $k$. Then every LinEMSOL($\tau _{1}$)
optimization problem on $mathcal{C}$ can be solved in linear time if a $k$%
-presentation of the graph is part of the input.
\end{theorem}

Liedloff et al.\cite{LKLP} used Theorem \ref{courcelle} to show that the
complexity of the Roman domination decision problem could be relaxed under
certain restrictions on the urdelying graphs. Now, we prove a similar result
regarding to the decision problem associated to the maximal double Roman
domination problem.

\begin{theorem}
\label{mdrdlinear} Problem MAXIMAL\ DOUBLE ROM-DOM belongs to the class of
optimization problems LinEMSOL($\tau _{1}$).
\end{theorem}

\noindent \textbf{Proof. } To prove the result we need to describe the
maximal double Roman domination problem as a LinEMSOL($\tau_1$) optimization
problem.

Let $f=(V_0, V_1,V_2,V_3)$ be a MDRD function with minimum weight in the
graph $G=(V,E)$. First, note that the weight of $f$ is $%
w(f)=|V_1|+2|V_2|+3|V_3|.$ Next, let us define the free set-variables $X_i:V
\rightarrow \{0,1\}$ as follows
\begin{equation*}
X_i(z)=1, \mbox{ if and only if } z \in V_i
\end{equation*}
and $X_i(z)=0$ in other case.

We use the notation $|X_i|=\sum_{z\in V} X_i(z)$ for consistency with
logical notations even if it is equivalent to $|X_i|=|V_i| .$

Taking into account the definition of a MDRD function, we can describe the
decision problem associated as

\begin{equation*}
\min_{X_i} \left\{ |X_1|+2|X_2|+3|X_3| \; : \; <G(\tau_1),X_0,\ldots, X_3 >
\; \vDash \theta(X_0,\ldots ,X_3) \right\},
\end{equation*}
where $\theta$ is defined below
\begin{equation*}
\begin{array}{c}
\theta(X_0,\ldots ,X_3) = \forall z \left( {\ \vphantom{\frac{\frac{A}{B}}{%
\frac{A}{B}}} } X_3(z) \vee X_2(z) \vee \left( {\vphantom{\frac{\int A}{B}}}
X_1 (z) \land \exists t \left( \left( X_2(t) \vee X_3(t) \right) \land
R(z,t) \right) {\vphantom{\frac{\int A}{B}}} \right) \vee \right. \\[0.5em]
\vee \left. \left( {\vphantom{\frac{\int A}{B}}} X_0(z) \land \left( \exists
t \left( X_3(t) \land R(z,y) \right) \vee \exists t,v \left( X_2(t) \land
X_2(v) \land R(z,t) \land R(z,v) \right) \right) \right) {\ %
\vphantom{\frac{\frac{A}{B}}{\frac{A}{B} }} } \right) \vee \\[0.5em]
\vee \; \exists z \left( {\ \vphantom{\frac{\frac{A}{B}}{\frac{A}{B}}} }
\left( {\vphantom{\frac{\int A}{B}}} X_1(z) \vee X_2(z) \vee X_3(z) {%
\vphantom{\frac{\int A}{B}}} \right) \land \forall t \left( {%
\vphantom{\frac{\int A}{B}}} R(z,t) \rightarrow X_1(t) \vee X_2(t) \vee
X_3(t) {\vphantom{\frac{\int A}{B}}} \right) {\ \vphantom{\frac{%
\frac{A}{B}}{\frac{A}{B}}} } \right)%
\end{array}%
\end{equation*}

Clearly, $\theta$ defines the corresponding propierties to a maximal double
roman domination funcion in the graph $G$. The definition of $\theta$ may be
divided into two main clauses, the first of which describes the conditions
for $f$ to be a double roman domination function. The second one assures
that there is at least a vertex with a possitive label that have no
neighbours in $V_0$. Therefore, $f$ is a maximal double roman domination
function in $G$ if and only if the logical expression $\theta$ is satisfied,
which finishes the prove. \hfil $\Box$

Next, as a direct consequence of the previous theorem, we can derive some
aditional results.

\begin{corollary}
Problem MAXIMAL\ DOUBLE ROM-DOM can be solved in linear time on any graph $G$
with clique-width bounded by a possitive intenger $k$, provided that either
there exists a linear-time algorithm to construct a $k$-expression of $G$,
or a $k$-expression of $G$ is part of the input.
\end{corollary}

Taking into account that any bounded treewidth graph is also a bounded
clique-width graph, we also deduce the following corollary.

\begin{corollary}
Problem MAXIMAL\ DOUBLE ROM-DOM can be solved in linear time for any tree
graph $G$.
\end{corollary}

{Other than trees that are well-known they have treewidth equal to 1, here
are two other families of graphs having bounded clique-width }$cw(G)${:
cographs }($cw(G)\leq 2$) and {}distance hereditary graphs ($cw(G)\leq 3$). {%
Therefore MAXIMAL\ DOUBLE\ ROM-DOM is also solvable in linear time }for all
these graphs.

\section{Relationships between $\protect\gamma _{dR}^{m}$ and $\protect%
\gamma _{dR},\protect\gamma _{mR},\protect\gamma _{m}$}

{In this section, we present some results relating the maximal DRD-number of
a graph to some domination parameters including the DRD-number, maximal
RD-number and the maximal domination number. }

\begin{proposition}
\label{=}For any connected graph $G$,
\begin{equation*}
\gamma _{dR}(G)\leq \gamma _{dR}^{m}(G),
\end{equation*}%
with equality if and only if $G$ satisfies one of the following.

\begin{enumerate}
\item $G$ is trivial.

\item $\delta (G)=1$ and there exists a $\gamma _{dR}(G)$-function that
assigns 1 to {a leaf and 2 to its support vertex}.

\item $\delta (G)=2$ and there exists a $\gamma _{dR}(G)$-function that
assigns 1 to a vertex of degree 2 and one of its neighbors.
\end{enumerate}
\end{proposition}

\textbf{Proof. } Since any $\gamma _{dR}^{m}(G)$-function is a {DRDF }of $G$%
, we have $\gamma _{dR}(G)\leq \gamma _{dR}^{m}(G),$

If $G$ is trivial, then clearly $\gamma _{dR}(G)=\gamma _{dR}^{m}(G)$. If $%
\delta (G)=1$ and there exists a $\gamma _{dR}(G)$-function $f$ which
assigns 1 to a {leaf }and 2 to its support vertex, then obviously $f$ is an MDRDF of $G
$ and this implies that $\gamma _{dR}(G)=\gamma _{dR}^{m}(G)$. {Likewise, if
}$\delta (G)=2$ and there exists a $\gamma _{dR}(G)$-function that assigns 1
to a vertex of degree 2 and one of its neighbors, {then clearly }$f$ {is an
MDRDF of }$G$ {and the result follows as above. }

Conversely, assume that $\gamma _{dR}(G)=\gamma _{dR}^{m}(G)$, {and let }$%
f=(V_{0},V_{1},V_{2},V_{3})$ be a $\gamma _{dR}^{m}(G)$-function. {Let }$%
v\in V_{1}\cup V_{2}\cup V_{3}$ be a vertex not dominated by $V_{0}$. If $%
f(v)=3$, then the function $g$ defined by $g(v)=2$ and $g(x)=f(x)$
otherwise, is a DRDF of $G$ of weight $\gamma _{dR}(G)-1$ which is a
contradiction. Hence $f(v)\in \{1,2\}$. We proceed with {the following }%
claim.

\smallskip \noindent \textbf{Claim.} $\deg (v)\leq 2$.\newline
\textbf{Proof of Claim.} Suppose, to the contrary, that $\deg (v)\geq 3$.
Since $v$ is not dominated by $V_{0}$, we have $N(v)\subseteq V_{1}\cup
V_{2}\cup V_{3}$. If $v$ has three neighbors with label 1 under $f$, then
the function $g$ defined by $g(v)=3$, $g(x)=0$ for $x\in \{y\in N(v)\mid
f(y)=1\}$ and $g(x)=f(x)$ otherwise, is a DRDF of $G$ of weight at most $%
\gamma _{dR}(G)-1,$ a contradiction. Hence $v$ has at most two neighbors
with label 1. First let $f(v)=2$. If $v$ has two neighbors neighbor with
label 1 under $f$, then the function $g$ defined before, is a DRDF of $G$ of
weight less than $\gamma _{dR}(G)$ which is a contradiction. Therefore $v$
has at most one neighbor with label 1, and {thus }$v$ has at least two
neighbors with label greater than 1. Then the function $h$ defined by $g(v)=0
$, $g(x)=2$ for $x\in \{y\in N(v)\mid f(y)=1\}$ and $g(x)=f(x)$ otherwise,
is a DRDF of $G$ of weight at most $\gamma _{dR}(G)-1,$ a contradiction. Now
let $f(v)=1$. If $v$ has two neighbors with label 1, then the function $h_{1}
$ defined by $h_{1}(v)=2$, $h_{1}(x)=0$ for $x\in \{y\in N(v)\mid f(y)=1\}$
and $h_{1}(x)=f(x)$ otherwise, is a DRDF of $G$ of weight most $\gamma
_{dR}(G)-1,$ a contradiction. Otherwise, $v$ has two neighbors with label
greater than 1 and the function $h_{2}$ defined by $h_{2}(v)=0$ and $%
h_{2}(x)=f(x)$, is a DRDF of $G$ of weight $\gamma _{dR}(G)-1,$ a
contradiction too. This proves the claim. $\hfill \blacklozenge $

Thus $\deg (v)\leq 2$. {Recall that } $f(v)\in \{1,2\},$ {and }consider {the
following two }cases.

\smallskip \noindent \textbf{Case 1.} $\deg (v)=1$.\newline
Let $u$ be the support vertex of $v$. If $f(v)=2$, then obviously $f(u)=1$ and the
function $g$ defined by $g(v)=1$ and $g(u)=2$ satisfies {item }(2). Hence
assume that $f(v)=1$. Since $f$ is {also a} $\gamma _{dR}(G)$-function, we
must have $f(u)=2$, and so $G$ {item }satisfies (2).

\smallskip \noindent \textbf{Case 2.} $\deg (v)=2$.\newline
Let $u,w$ be the neighbors of $v$. If $f(u),f(w)\geq 2$, then the function $%
h_{2}$ defined {in the proof of the claim, }is a DRDF of $G$ of weight {at
most }$\gamma _{dR}(G)-1$, a contradiction. {Hence }we may assume, without
loss of generality, that $f(u)=1$. We claim that $f(v)=1$. Suppose, to the
contrary, that $f(v)=2$. If $f(w)=1$, then the function $g$ defined by $%
g(v)=3$, $g(u)=g(w)=0$ and $g(x)=f(x)$ otherwise, is a is a DRDF of $G$ of
weight $\gamma _{dR}(G)-1,$ a contradiction. If $f(w)\geq 2$, then the
function $g$ defined by $g(v)=0$, $g(u)=2$ and $g(x)=f(x)$ otherwise, is a
DRDF of $G$ of weight $\gamma _{dR}(G)-1,$ a contradiction again. Thus $%
f(v)=1$ {and since }$f(u)=1$ {item (3) follows. } $\ \Box $

\begin{proposition}
\label{mR1}For any graph $G$, $\gamma _{dR}^{m}(G)\leq 2\gamma _{mR}(G)$
with equality if and only if $G=\overline{K_{n}}$.
\end{proposition}

\textbf{Proof. } Let $f=(V_{0},V_{1},V_{2})$ be a $\gamma _{mR}(G)$-function
that minimizes the number of vertices in $V_{1}$. {Note that }$\gamma
_{mR}(G)=|V_{1}|+2|V_{2}$. {Clearly, }$(V_{0},\emptyset ,V_{1},V_{2})$ {is
an MDRDF of }$G,$ {and thus }$\gamma _{dR}^{m}(G)\leq
2|V_{1}|+3|V_{2}|=\gamma _{mR}(G)+|V_{1}|+|V_{2}|\leq 2\gamma _{mR}(G)$. {%
Now, if }$\gamma _{dR}^{m}(G)=2\gamma _{mR}(G)$, then we must have {equality
throughout the previous inequality chain, and thus }$V_{2}=\emptyset $.
Hence, $V_{0}=\emptyset $ must hold, and so $V=V_{1}$. Since $|V_{1}|$ is
minimized under $f$ , we deduce that {each component of }$G$ {has order at
most two. Now, assume that }$G$ {has a component of order two and let }$u$
and $v$ {be the vertices of such a component. Then }function $g$ which
assigns 1 to $u$, 2 to $v$, and 2 to every other vertex is a MDRDF of $G$ {%
of weight }$2n-1<2\gamma _{mR}(G)$, a contradiction. {Therefore each
component of }$G$ {is trivial and thus }$G=\overline{K_{n}}$. $\ \Box $

\bigskip

Restricted to isolated-free graphs $G${, it follows from
Proposition \ref{mR1} that }$\gamma _{dR}^{m}(G)\leq 2\gamma _{mR}(G)-1.$ In
the next we characterize isolated-free graphs attaining this upper bound. {\
}

\begin{corollary}
\label{mR}If $G$ is an isolated-free graph, then $\gamma _{dR}^{m}(G)\leq
2\gamma _{mR}(G)-1,$ with equality if and only if {\ $G=K_{2}.$}
\end{corollary}

\textbf{Proof. } Assume that $\gamma _{dR}^{m}(G)=2\gamma _{mR}(G)-1.$ Let
the components of $G$ be $G_{1},G_{2},...,G_{p},$ and let $%
f=(V_{0},V_{1},V_{2})$ be a $\gamma _{mR}(G)$-function such that $\left\vert
V_{1}\right\vert $ is as small as possible. Recall that $\gamma
_{mR}(G)=|V_{1}|+2|V_{2}|.$ Now, since $\gamma
_{dR}^{m}(G)=\sum_{i=1}^{p}\gamma _{dR}^{m}(G_{i}),$ $\gamma
_{mR}(G)=\sum_{i=1}^{p}\gamma _{mR}(G_{i})$ and for each $i$, $\gamma
_{dR}^{m}(G_{i})\leq 2\gamma _{mR}(G_{i})-1,$ we deduce that $p=1,$ that is $%
G$ is connected. Moreover, since $(V_{0},\emptyset ,V_{1},V_{2})$ is an
MDRDF of $G,$ we have

\begin{equation*}
2\gamma _{mR}(G)-1=\gamma _{dR}^{m}(G)\leq 2|V_{1}|+3|V_{2}|=2\gamma
_{mR}(G)-|V_{2}|,
\end{equation*}%
and thus $|V_{2}|\leq 1$. First, assume that $|V_{2}|=0.$ Then $%
V_{0}=\emptyset $, and clearly $V_{1}=V(G).$ Hence $\gamma _{mR}(G)=n.$ We
claim that $G$ is nontrivial complete graph. Suppose not and let $u$ and $v$
be two non adjacent vertices of $G.$ Then assigning $0$ to $u,$ $2$ to a
neighbor of $u,$ and $1$ to the remaining vertices of $G$ provides an MRDF
of $G$ with less vertices { assigned} $1$ than under $f,$ a
contradiction. Hence $G=K_{n},$ with $n\geq 2.$ {Now, since }$\gamma
_{dR}^{m}(K_{n})=n+1$ and {$\gamma _{dR}^{m}(G)=2\gamma _{mR}(G)-1$ we
deduce from $n+1=2n-1$ that $G=K_{2}$.}

Assume now that $|V_{2}|=1.$ Then $\gamma _{mR}(G)=|V_{1}|+2$ and thus $%
\gamma _{dR}^{m}(G)=2|V_{1}|+3.$ Denote by $V_{1}^{\prime }=\{w\in
V_{1}:N(w)\cap V_{2}\neq \emptyset \}$ and $V_{1}^{\prime \prime
}=V_{1}\!\smallsetminus \!V_{1}^{\prime }.$ Clearly, since $f$ is an MRDF of
$G,$ there is a vertex, say $t\in V_{1}$ such that $N(t)\cap V_{0}=\emptyset
.$ Therefore the function $g=(V_{0},V_{1}^{\prime },V_{1}^{\prime \prime
},V_{2})$ is an MDRDF of $G$ and thus
\begin{equation*}
2|V_{1}|+3=\gamma _{dR}^{m}(G)\leq w(g)=|V_{1}^{\prime }|+2|V_{1}^{\prime
\prime }|+3=|V_{1}|+|V_{1}^{\prime \prime }|+3
\end{equation*}%
Therefore, $V_{1}=V_{1}^{\prime \prime }$ and hence $V_{1}^{\prime
}=\emptyset .$ Now, since $G$ is connected, we must have $V_{1}\cap N(t)\neq
\emptyset .$ But then the function $h=(V_{0},\{t\},V_{1}-\{t\},V_{2})$ is an
MDRDF of $G$ and thus
\begin{equation*}
2|V_{1}|+3=\gamma _{dR}^{m}(G)\leq w(h)=1+2(|V_{1}|-1)+3=2|V_{1}|+2,
\end{equation*}%
which is a contradiction.

The converse is obvious. \hfill $\Box $

\bigskip

As consequence of Proposition {\ref{mR1} and Corollary }\ref{mR}, the
following is immediate.

\begin{corollary}
If $G$ is a nontrivial connected graph different from $K_{n},$ then $\gamma
_{dR}^{m}(G)\leq 2\gamma _{mR}(G)-2.$
\end{corollary}

Next we show that {for every graph }$G,$ the maximal {RD-number }is strictly
smaller than the maximal {DRD-number.}

\begin{proposition}
\label{mR2}For every graph $G$, $\gamma _{mR}(G)<\gamma _{dR}^{m}(G)$.
\end{proposition}

\textbf{Proof. } Let $f=(V_{0},V_{1},V_{2},V_{3})$ be any $\gamma
_{dR}^{m}(G)$-function. If $V_{3}\neq \emptyset $, then every vertex in $%
V_{3}$ can be reassigned the value 2 and the resulting function will be {an
MRDF of }$G,$ {implying that }$\gamma _{mR}(G)<\gamma _{dR}^{m}(G)$. {Hence a%
}ssume\ that $V_{3}=\emptyset $. By definition $V_{2}\neq \emptyset $. {Now,
if }$V_{0}=\emptyset $, then every vertex in $V_{2}$ can be reassigned the
value 1 and the resulting function will be {an MRDFof }$G,$ {implying again
that }$\gamma _{mR}(G)<\gamma _{dR}^{m}(G)$. {Thus a}ssume that $V_{0}\neq
\emptyset ,$ {and let }$v\in V_{2}$. Since every vertex in $V_{0}$ is
adjacent to at least two vertices in $V_{2}$, then by reassigning $v$ {the
value }1 {provides an MRDF }on $G$, {which implies that }$\gamma
_{mR}(G)<\gamma _{dR}^{m}(G)$. $\ \Box $

\bigskip

{\ According to Corollary \ref{mR} and Proposition\ref{mR2}, we have: }

\begin{corollary}
For any nontrivial connected graph $G$, $\gamma _{mR}(G)<\gamma
_{dR}^{m}(G)<2\gamma _{mR}(G)$.
\end{corollary}

Since $V_{1}\cup V_{2}\cup V_{3}$ is a maximal dominating set when $%
f=(V_{0},V_{1},V_{2},V_{3})$ is an MDRDF, and since placing a 3 at {each
vertex }of a maximal dominating set {and a 0 elsewhere} yields an MDRDF, we
obtain the following result.

\begin{obs}
\label{eqq} For any graph $G$, \
\begin{equation*}
\gamma _{m}(G)\leq \gamma _{dR}^{m}(G)\leq 3\gamma _{m}(G).
\end{equation*}
\end{obs}

{\ Our next result slightly improves the upper bound of Observation \ref{eqq}%
. }

\begin{proposition}
\label{prop_gm}For any graph $G$,
\begin{equation*}
\gamma _{m}(G)+1\leq \gamma _{dR}^{m}(G)\leq 3\gamma _{m}(G)-1.
\end{equation*}%
Furthermore, these bounds are sharp.
\end{proposition}

\textbf{Proof. } {We first }prove the lower bound. Let $%
f=(V_{0},V_{1},V_{2},V_{3})$ be a $\gamma _{dR}^{m}$-function of $G$.
Clearly, $V_{2}\cup V_{3}\neq \emptyset $ and $V_{1}\cup V_{2}\cup V_{3}$ {%
is a maximal dominating set of }$G.$ {Hence, }$\gamma
_{dR}^{m}(G)=|V_{1}|+2|V_{2}|+3|V_{3}|\geq |V_{1}|+|V_{2}|+|V_{3}|+1\geq
\gamma _{m}(G)+1$. This bound is sharp complete graphs.

{Now} we prove the upper bound. Let $S$ be a maximal dominating set of $G$.
If $V=S$, then the function $f$ defined on $G$ by $f(x)=2$ for each $x\in V$%
, is {an} MDRDF of $G$ and so $\gamma _{dR}^{m}(G)\leq 2\gamma _{m}(G)\leq
3\gamma _{m}(G)-1,$ as desired. Hence we assume that $V-S\neq \emptyset $.
Since $V-S$ is not a dominating set {of }$G$, there is a vertex $v\in S$
which is not dominated by $V-S$. Then the function $f$ defined on $G$ by $%
f(v)=2$, $f(x)=3$ for each $x\in S-\{v\}$ and $f(x)=0$ otherwise, is {an}
MDRDF of $G$ and {thus }$\gamma _{dR}^{m}(G)\leq 3\gamma _{m}(G)-1.$

To {see }the sharpness {of the upper bound, }let $G$ be {a graph }obtained
from a cycle $C_{n}=(v_{1}v_{2}\dots v_{n}v_{1})$ by {first adding for each }%
$v_{i}$, $t\geq 2$ {new vertices }$v_{i}^{1},\ldots ,v_{i}^{t}$ {attached by
}edges $v_{i}v_{i}^{1},\ldots ,v_{i}v_{i}^{t}$, {and then adding an }%
isolated vertex $v$. {One can easily see }that $S=\{v,v_{1},\dots ,v_{n}\}$
is a minimum maximal dominating set of $G$ and so $\gamma _{m}(G)=n+1$. {%
Moreover, }one can {also }see that the function $f$ defined by $f(v)=2$, $%
f(v_{1})=\dots =f(v_{n})=3$ and $f(x)=0$ otherwise is the unique $\gamma
_{dR}^{m}(G)$-function and so $\gamma _{dR}^{m}(G)=3n+2=3\gamma _{m}(G)-1.$ $%
\ \Box $

\bigskip

{The next result slightly improves the upper bound of Proposition \ref%
{prop_gm} for isolated-free graphs. }

\begin{proposition}
For any isolated-free graph $G$,
\begin{equation*}
\gamma_{dR}^m(G)\le 3\gamma_m(G)-2.
\end{equation*}
Furthermore, this bound is sharp.
\end{proposition}

\textbf{Proof. }Let $S$ be a maximal dominating set of $G$. {Note that since
}$G$ {has no isolated vertices, }$\left\vert S\right\vert \geq 2.$ If $V=S$,
then the function $f$ defined on $G$ by $f(x)=2$ for each $x\in V$, is a
MDRDF of $G$ and so $\gamma _{dR}^{m}(G)\leq 2\gamma _{m}(G)\leq 3\gamma
_{m}(G)-2$. {Hence let }$V-S\neq \emptyset $. Since $V-S$ is not a
dominating set {of }$G$, {let }$v\in S$ {be a vertex }not dominated by $V-S$%
. {Clearly, }$v$ {has t least one neighbor in }$S$ {because of }$G$ {is
isolated-free. }Then the function $f$ defined on $G$ by $f(v)=1$, $f(x)=3$
for each $x\in S-\{v\}$ and $f(x)=0$ otherwise, is a MDRDF of $G$ and so $%
\gamma _{dR}^{m}(G)\leq 3\gamma _{m}(G)-2.$

To show the sharpness {of the bound}, let $G$ be {a graph }obtained from a
cycle $C_{n}=(v_{1}v_{2}\dots v_{n}v_{1})$ by {first adding for each }$v_{i}$%
, $t\geq 3$ {new vertices }$v_{i}^{1},\ldots ,v_{i}^{t}$ {attached by }edges
$v_{i}v_{i}^{1},\ldots ,v_{i}v_{i}^{t}$. It is not hard to see that the set $%
S=\{v_{1}^{1},v_{1},\dots ,v_{n}\}$ is a minimum maximal dominating set of $%
G $ and so $\gamma _{m}(G)=n+1$. {Moreover, }the function $f$ defined by $%
f(v_{1}^{1})=1$, $f(v_{1})=\dots =f(v_{n})=3$ and $f(x)=0$ otherwise is a $%
\gamma _{dR}^{m}(G)$-function and so $\gamma _{dR}^{m}(G)=3n+1=3\gamma
_{m}(G)-2.$ $\ \Box $

\begin{proposition}
\label{delta}For any graph $G$ without isolated vertices,
\begin{equation*}
\gamma _{dR}^{m}(G)\leq \gamma _{dR}(G)+\delta (G).
\end{equation*}%
Furthermore, this bound is sharp.
\end{proposition}

\textbf{Proof.} Let $f=(V_{0},V_{1},V_{2},V_{3})$ be a $\gamma _{dR}(G)$%
-function. If $V_{0}=\emptyset $ or $V_{0}$ does not dominate $V_{1}\cup
V_{2}\cup V_{3},$ then $f$ is an MDRDF of $G,$ and the result is clearly
valid. Hence we assume that $V_{0}$ is non-empty and dominates all $V(G).$
Let $v$ be a vertex of minimum degree in $G$ and let $A_{v}=N[v]\cap V_{0}.$
Since $V_{0}$ is a dominating set of $G,$ $\left\vert A_{v}\right\vert \leq
\delta (G).$ It follows that the function $g=(V_{0}-A_{v},V_{1}\cup
A_{v},V_{2},V_{3})$ is an MDRDF on $G,$ and thus $\gamma _{dR}^{m}(G)\leq
\omega (g)=\omega (f)+\left\vert A_{v}\right\vert \leq \gamma
_{dR}(G)+\delta (G)$.

To see the sharpness of the upper bound, let $G$ be the graph obtained from $%
K_{n}$ ($n\geq 2$) by adding a new vertex attached by an edge to exactly one
vertex of $K_{n}$. Then $\gamma _{dR}(G)=3$ and $\gamma _{dR}^{m}(G)=\gamma
_{dR}(G)+1$. $\ \Box $

\begin{proposition}
Let $G$ be a connected graph of order $n$ with $\mathrm{diam}(G)\geq 4$.
Then
\begin{equation*}
\gamma _{dR}^{m}(G)\leq 2(n-\delta (G)).
\end{equation*}
\end{proposition}

\textbf{Proof.} Let $P=u_{1}u_{2}\ldots u_{\mathrm{diam}(G)+1}$ be a
diametral path in $G$, {and consider the }function $f=(N(u_{2}),\{u_{5}%
\},V(G)-(N[u_{2}]\cup \{u_{5}\}),\{u_{2}\}).$ {Then }$f$ is an MDRDF of $G$
and thus
\begin{equation*}
\begin{array}{clc}
\gamma _{dR}^{m}(G) & \leq & \omega (f)\hfill \\
& = & |V_{1}|+2|V_{2}|+3|V_{3}|\hfill \\
& = & 2(n-\deg (u_{2})-2)+3+1\hfill \\
& \leq & 2n-2\delta (G),\hfill%
\end{array}%
\end{equation*}%
and the proof is complete. $\ \Box $

\section{Maximal DRD-number in special graphs}

{In this section, we determine the exact values of the maximal DRD-number
for paths and cycles as well as an upper bound for trees in terms of the
order.\ We begin by the following corollary which is immediate from
Proposition \ref{=} and \ref{delta}. }

\begin{corollary}
\label{cor1}For every nontrivial tree $T,$
\begin{equation*}
\gamma _{dR}(G)\leq \gamma _{dR}^{m}(G)\leq \gamma _{dR}(G)+1.
\end{equation*}
\end{corollary}

%
{Let $\mathcal{F}$ be the family of all trees that can be built from {$k\geq
1$ paths $P_{4}^{i}:=v_{1}^{i}v_{2}^{i}v_{3}^{i}v_{4}^{i}\;(1\leq i\leq k)$
by adding $k-1$ edges incident with the $v_{2}^{i}$'s so that they induce a
connected subgraph.} In \cite{bhh}, Beeler et al. gave an upper bound on the
DR-number for trees in terms of their order and characterized the trees
reaching this bound. }

\begin{theorem}
\label{tree1}If $T$ is a tree with order $n\geq 3$, then $\gamma
_{dR}(T)\leq 5n/4$ with equality if and only if $T\in \mathcal{F}$.
\end{theorem}

{Clearly, by Corollary \ref{cor1} and Theorem \ref{tree1}, every tree }$T$ {%
of order }$n\geq 3$, {satisfies }$\gamma _{dR}^{m}(G)\leq 5n/4+1.$ {However,
we will see with the next result that the }$\frac{4}{5}n$ {by upper bound
remains valid for the maximal DRD-number of trees of order at least four.
Recall that a tree is \emph{double star }if it contains exactly vertices
that are not leaves. Moreover, a double star with respectively }$p$ and $q$ {%
leaves attached at each support vertex is denote by }$S_{p,q}.$

\begin{theorem}
\label{tree}If $T$ is a tree with order $n\geq 4$, then $\gamma
_{dR}^{m}(T)\leq 5n/4$ with equality if and only if $T\in \mathcal{F}.$
\end{theorem}

\textbf{Proof. } Let $T$ be a tree with order $n\geq 4$. We will proceed by
induction on $n$. Since $n\geq 4$, $\mathrm{diam}(T)\geq 2$. If the diameter
of $T$ is $2$, then $T$ is the star $K_{1,n-1}$ for $n\geq 4$ and we have $%
\gamma _{dR}^{m}(T)=4<5n/4$. {Hence assume that }$T$ {has diameter }$3,$ {%
that is, }$T$ is a double star $S_{r,s}$ for $1\leq r\leq s$. {Clearly, }$%
n=r+s+2\geq 4$. {Now, if }$s=1$, then $T=P_{4},$ {where }$\gamma
_{dR}^{m}(T)=5=5n/4$ {and }$T\in \mathcal{F}$. {If }$s\geq 2$, {then }$%
\gamma _{dR}^{m}(T)=6<5n/4$ {when }$r=1$ {while }$\gamma _{dR}^{m}(T)=7<5n/4$
{when }$r\geq 2.$ {Hence we }may suppose that the diameter of $T$ is at
least 4. This implies that $n\geq 5$. {Moreover, let }the result hold for
any tree $T^{\prime }$ with order $4\leq n^{\prime }<n$. Among all longest
paths in $T$, choose $P=v_{1}\ldots ,v_{k}$ to be one that maximizes the
degree of $v_{2}$. {Note that }$k\geq 5$ {and }$\deg (v_{k-1})\leq \deg
(v_{2}).$ Root T at $v_{k}$. Note that by our choice, every child of $v_{2}$
is a leaf. We consider three cases.

\smallskip \noindent \textbf{Case 1.} $\deg (v_{2})>3$.\newline
Then $v_{2}$ {is a support vertex with at least three leaves}. {Let }$T^{\prime
}=T-v_{1}$. {Clearly, there is a }$\gamma _{dR}(T^{\prime })$-function $f$
on $T^{\prime }$ {that assigns }$3$ to $v_{2}$. Hence, assigning $v_{1}$ a
1, the {function} $f$ {can be extended to be an MDRDF }of $T$ , implying
that $\gamma _{dR}^{m}(T)\leq \gamma _{dR}(T^{\prime })+1\leq 5(n-1)/4+1<5n/4
$.

\smallskip \noindent \textbf{Case 2.} $\deg (v_{2})=3$.\newline
Let $v_{1},w$ be the leaf neighbors of $v_{3}$ and {let }$T^{\prime
}=T-T_{v_{2}}$. Since $\mathrm{diam}(T)\geq 4$, we have $n(T^{\prime })\geq 3
$. {Now if }$n(T^{\prime })=3$, then the function $g$ defined on $T$ by $%
g(v_{2})=3,g(v_{4})=2,g(v_{5})=1$ and $g(x)=0$ otherwise, is an MDRDF of $T$
and so $\gamma _{dR}^{m}(T)\leq 6<5n/4$. {Hence assume }that $n(T^{\prime
})\geq 4$, {and let }$f$ be a $\gamma _{dR}^{m}$-function of $T^{\prime }$. {%
Then }assigning $v_{2}$ a 3, and $v_{1},w$ a 0, the {function} $f$ {can be
extended to be an MDRDF }of $T$ , implying that $\gamma _{dR}^{m}(T)\leq
\gamma _{dR}^{m}(T^{\prime })+3\leq 5(n-3)/4+3<5n/4$.

\smallskip \noindent \textbf{Case 3.} $\deg (v_{2})=2$.\newline
{Clearly by our choice of }$P,$ $\deg (v_{k-1})=\deg (v_{2}).$ We consider
two subcases.

\smallskip \textbf{Subcase 3.1.} $\deg (v_{3})=2$.\newline
First let $\deg (v_{4})\geq 3$. Let $T^{\prime }=T-T_{v_{3}}$. {Observe that
since }$\mathrm{diam}(T)\geq 4$ {and} $\deg (v_{k-1})=\deg (v_{2}),$ {we
have }$n(T^{\prime })\geq 4$. Applying our inductive hypothesis, we have $%
\gamma _{dR}^{m}(T^{\prime })\leq 5(n-3)/4$. {Now, let }$f$ be a $\gamma
_{dR}^{m}$-function. If $f(v_{4})=0$ or $f(v_{4})\geq 1$ and $v_{4}$ is
dominated by $V_{0}^{f}$, then by assigning 3 to $v_{2}$, and a 0 to each of
$v_{1}$ and $v_{3}$, {function} $f$ {can be extended to be an MDRDF }of $T$,
implying that $\gamma _{dR}^{m}(T)\leq \gamma _{dR}^{m}(T^{\prime })+3\leq
5(n-3)/4+3<5n/4$. {Hence let }$f(v_{4})\geq 1$ and $v_{4}$ is not dominated
by $V_{0}$. If $f(v_{4})\geq 2$, then by assigning 2 to $v_{2}$, 1 to $v_{1}$
and 0 to $v_{3}$, we {obtain an MDRDF }of $T$, {and clearly }$\gamma
_{dR}^{m}(T)<5n/4$. {Now let }$f(v_{4})=1$. If $v_{4}$ has a neighbor with
label 3 or two neighbors with label 3, then by {reassigning }$v_{4}$ {a 0,
and }assigning 3 to $v_{2}$, 1 to $v_{1}$ and 0 to $v_{3}$, we {obtain an
MDRDF }of $T$, implying that $\gamma _{dR}^{m}(T)<5n/4$. {Hence a}ssume that
$v_{4}$ has exactly one neighbor with label 2. {Since }$\deg (v_{4})\geq 3$,
{vertex }$v_{4}$ has at least one neighbor with label 1, say $w$. Then by {%
reassigning }$v_{4}$ and $w$ {the values 2 and 0, respectively, and
assigning }2 to $v_{2}$, 1 to $v_{1}$ and 0 to $v_{3}$, we {obtain an MDRDF }%
of $T$ {of weight less than }$5n/4$.

\smallskip \noindent Now let $\deg (v_{4})=2$. Let $T^{\prime }=T-T_{v_{4}}$%
, {and consider a }$\gamma _{dR}(T^{\prime })$-function $f$. {If }$T^{\prime
}$ {has order 1, then }$T=P_{5}$ {and }$\gamma _{dR}^{m}(T)=6<25/4=5n/4$. If
$T^{\prime }$ has order 2, then $T=P_{6}$ and $\gamma
_{dR}^{m}(T)=7<30/4=5n/4$. If $T^{\prime }$ has order 3, then {since }$\deg
(v_{k-1})=\deg (v_{2}),$ {we have }$T=P_{6}$ {and so one can see that }$%
\gamma _{dR}^{m}(T)\leq 8<35/4=5n/4$. {Therefore, }we may assume that $%
n^{\prime }\geq 4$. If $T\not\in \mathcal{F}$, then by Theorem \ref{tree1},
we have $\gamma _{dR}(T)<5(n-4)/4$, and by assigning 2 to $v_{2},v_{4}$, 1
to $v_{1}$ and 0 to $v_{3}$, we {obtain an MDRDF }of $T$ {of weight less
than }$5n/4$. Thus, assume that $T^{\prime }\in \mathcal{F}$. Then we may
assume that $f(v_{5})\geq 2$. By by assigning {\ 3 }to $v_{2},$ 1 to $v_{4}$
and \ 0 to each of $\ v_{1},v_{3}$, we {obtain an MDRDF }of $T$ {of weight
less than }$5n/4$.

\smallskip \textbf{Subcase 3.2.} $\deg (v_{3})\geq 3$.\newline
By our choice of diametral path, every child of $v_{3}$ in $T_{v_{3}}$ is
either a leaf or a {support vertex of degree two}. Let $T^{\prime }=T-T_{v_{3}}$. If $%
T^{\prime }$ has order 2, then $T$ is a {tree obtained from a star }$%
K_{1,\deg (v_{3})}$ {by subdividing at least two edges of the star. In this
case, one can easily see that }$\gamma _{dR}^{m}(T)<5n/4$. If $n(T^{\prime
})=3$, then {one can also see }that $\gamma _{dR}^{m}(T)<5n/4$. So, we may
assume that $n^{\prime }\geq 4$. Let $r$ {and }$s$ be the numbers of
children of $v_{3}$ {which are leaves and support vertices, respectively. Note that }$%
s\geq 1$ {because of }$v_{2}.${\ }First let $r+s\geq 3$ or $r=0$. Clearly
any $\gamma _{dR}^{m}$-function of $T^{\prime }$ can be extended to an MDRDF
of $T$ by assigning 2 to $v_{3}$ and {every }leaf of $T_{v_{3}}$ at distance
2 from $v_{3}$, 1 to every leaf neighbor of $v_{3}$ {(if any) }and 0 to
other children of $v_{3}$. {Therefore, }$\gamma _{dR}^{m}(T)<5n/4$. {Now we
can assume }that $r=s=1.$ {Let }$w$ be the leaf neighbor of $v_{3}$. If $%
T\not\in \mathcal{F}$, then any $\gamma _{dR}$-function of $T^{\prime }$ can
be extended to an MDRDF of $T$ by assigning 2 to $v_{3},v_{1}$, 1 to $w$ and
0 to $v_{2}$. It follows from Theorem \ref{tree1} that $\gamma
_{dR}^{m}(T)\leq \gamma _{dR}(T^{\prime })+5<5(n-4)/4+5<5n/4$. Assume that $%
T\in \mathcal{F}$. Then $T^{\prime }$ can be built from $k\geq 1$ paths $%
P_{4}^{i}:=v_{1}^{i}v_{2}^{i}v_{3}^{i}v_{4}^{i}\;(1\leq i\leq k)$ by adding $%
k-1$ edges incident with the $v_{2}^{i}$'s so that they induce a connected
subgraph. We may assume, without loss of generality, that $v_{4}\in
V(P_{4}^{k})$. If $v_{4}$ is a leaf of $T^{\prime }$, then $T^{\prime }$ has
a $\gamma _{dR}(T^{\prime })$-function that assigns 1 to $v_{5}$, and the
function $g$ defined by $g(v_{3})=g(v_{1})=2$, $g(v_{4})=g(v_{2})=0$, $%
g(w)=1 $ and $g(x)=f(x)$ otherwise, is an MDRDF of $T$ {of weight less than }%
$5n/4$. Hence $v_{4}$ is a support vertex. If $k=1$, then clearly $T\in \mathcal{F}$
and we have $\gamma _{dR}^{m}(T)=10=5n/4$. {Thus, a}ssume that $k\geq 2$. If
$v_{4}=v_{3}^{k}$, then the function $g$ defined by $g(v_{2}^{i})=3$ for $%
1\leq i\leq k-1$, $g(v_{4}^{i})=2$ for $1\leq i\leq k$, $%
g(v_{1}^{k})=g(v_{1})=g(v_{3})=2$, $g(w)=1$ and $g(x)=0$ otherwise, is an
MDRDF of $T$ {of weight less than }$5n/4$. Thus $v_{4}=v_{2}^{k}$ and so $%
T\in \mathcal{F}$. {In this case, }the function $g$ defined by $%
g(v_{2}^{i})=3$ for $1\leq i\leq k$, $g(v_{4}^{i})=2 $ for $1\leq i\leq k$, $%
g(w)=g(v_{2})=2$, $g(v_{1})=1$ and $g(x)=0$ otherwise, is an MDRDF of $T$ {%
of weight }$5n/4=\gamma _{dR}(T)\leq \gamma _{dR}^{m}(T)\leq 5n/4$.
Therefore $\gamma _{dR}^{m}(T)=5n/4$ and the proof is complete. $\ \Box $

\bigskip


{In the aim to establish the exact values of the Maximal DRD-number for
paths and cycles, we need the use of the following }result {that can be
found }in \cite{acs1}.

\begin{proposition}
\label{path} For $n\geq 1$, $\gamma _{dR}(P_{n})=\left\{
\begin{array}{lll}
n & \mathrm{if} & n\equiv 0\;(\mathrm{mod}\;3) \\
n+1 & \mathrm{if} & n\equiv 1\;\mathrm{or}\;2\;(\mathrm{mod}\;3).%
\end{array}%
\right. $
\end{proposition}

\begin{proposition}
\label{path1} For $n\geq 1$, $\gamma _{dR}^{m}(P_{n})=n+1$.
\end{proposition}

\textbf{Proof. } The result is trivial for $n\in \{1,2,3\}.$ {Thus, assume
that }$n\geq 4$, and let $P_{n}=v_{1}v_{2}\dots v_{n}$ be a path on $n$
vertices. First let $n\equiv 1\pmod 3$. If $n$ is odd, then the function $f$
defined by $f(v_{1})=f(v_{n})=1$, $f(v_{2i})=2$ for $1\leq i\leq (n-1)/2$
and $f(x)=0$ otherwise, is an MDRDF of $P_{n}$ of weight $n+1. $ If $n$ is
even, then the function $f$ defined by $f(v_{1})=1$, $f(v_{2i})=2 $ for $%
1\leq i\leq n/2$ and $f(x)=0$ otherwise, is an MDRDF of $P_{n}$ of weight $%
n+1.$ {In either case, }$\gamma _{dR}^{m}(P_{n})\leq n+1$. It follows from
Proposition \ref{path} that $n+1=\gamma _{dR}(P_{n})\leq \gamma
_{dR}^{m}(P_{n})\leq n+1$ and {the equality follows. }

Now let $n\equiv 2\pmod 3$. If $n$ is even, then the function $f$ defined by
$f(v_{1})=1$, $f(v_{2i})=2$ for $1\leq i\leq n/2$ and $f(x)=0$ otherwise, is
an MDRDF of $P_{n}$ of weight $n+1$. If $n$ is odd, then the function $f$
defined by $f(v_{1})=f(v_{n})=1$, $f(v_{2i})=2$ for $1\leq i\leq (n-1)/2$
and $f(x)=0$ otherwise, is an MDRDF of $P_{n}$ of weight $n+1.$ {As before,
by using Proposition \ref{path}, }we get $\gamma _{dR}^{m}(P_{n})=n+1$.

Finally, let $n\equiv 0\pmod 3$. {Clearly, }the function $f$ defined by $%
f(v_{3i+2})=3$ for $0\leq i\leq (n-3)/3$ and $f(x)=0$ otherwise, is the
unique $\gamma _{dR}(P_{n})$-function. {Since every }$\gamma _{dR}^{m}(P_{n})
$-function $g$ {is a DRDF }of $P_{n}$ {but not a }$\gamma _{dR}(P_{n})$%
-function, {we deduce that }$\gamma _{dR}^{m}(P_{n})=\omega (g)>\gamma
_{dR}^{m}(P_{n})+1=n+1$. {The equality follows from the following MDRDF }$g$
of weight $n+1$ defined by $g(v_{n})=1$, $g(v_{3i+2})=3$ for $0\leq i\leq
(n-3)/3$ and $f(x)=0$ otherwise. $\ \Box $

\begin{proposition}
\label{cycle1} For $n\geq 3$, $\gamma _{dR}^{m}(C_{n})=n+1.$
\end{proposition}

\textbf{Proof. } Let $C_{n}=(v_{1}v_{2}\dots v_{n}v_{1})$ be a cycle on $n$
vertices {and let }$f=(V_{0},V_{1},V_{2},V_{3})$ be an MDRDF of $C_{n}$ such
that {the total weight of vertices not }dominated by $V_{0}$ {is as small as
possible}. Without loss of generality, {assume that }$v_{1}$ is not
dominated by $V_{0}$. {It follows that }$f(v_{1}),f(v_{2}),f(v_{n})\geq 1$ {%
and }$f(v_{1})\leq 2$. If $f(v_{1})=1$, then we may assume that $f(v_{2})=2.$
{Since }the function $f$ {remains an MDRDF of the path }$P_{n}$ {resulting
from }$C_{n}$ {by the deletion of the edge }$v_{1}v_{n},$ {we deduce from }%
Proposition \ref{path1} {that }$\gamma _{dR}^{m}(C_{n})=\omega (f)\geq
\gamma _{dR}^{m}(P_{n})=n+1$. {Assume now that }$f(v_{1})=2$. {By our }%
choice of $f,$ we must have $f(v_{n})=f(v_{1})=1$. {Consider the path }$%
P_{n-3}$ {obtained from }$C_{n}$ {by the deletion of vertices }$%
v_{1},v_{2},v_{n}.$ {Clearly, the restriction of }$f$ {on }$V(P_{n-3})$ is
an MDRDF {and thus by }Proposition \ref{path1}, we have $\gamma
_{dR}^{m}(C_{n})=\omega (f)\geq \gamma _{dR}^{m}(P_{n-3})+3=n+2$. {In either
case, we have }$\gamma _{dR}^{m}(C_{n})\geq n+1$.

Next we prove the inverse inequality. If $n\equiv 1,3,5\pmod 6$, then the
function $g$ defined by $g(v_{1})=g(v_{n})=1$, $g(v_{2i})=2$ for $1\leq
i\leq (n-1)/2$ and $g(x)=0$ otherwise, is an MDRDF of $C_{n}$ of weight $n+1$%
. If $n\equiv 0,2,4\pmod 6$, then the function $g$ defined by $g(v_{1})=1$, $%
g(v_{2i})=2$ for $1\leq i\leq n/2$ and $g(x)=0$ otherwise, is an MDRDF of $%
C_{n}$ of weight $n+1$. {In either case,} $\gamma _{dR}^{m}(C_{n})\geq n+1$,
{and the desired equality follows. } $\ \Box $\\

\textbf{Acknowledgements}\newline

H. Abdollahzadeh Ahangar was supported by the Babol Noshirvani University of
Technology under research grant number BNUT/385001/00.


\end{document}